\documentclass[12pt,amsbsy]{article}
\usepackage{amsmath,amssymb,mathrsfs}
\newcommand{\xr}{X_1,\dots ,X_n}
\newcommand{\lejn}{1\leq j\leq n}
\newcommand{\rn}{random }
\newcommand{\fr}{free independence }
\newcommand{\alg}{algebra}
\newcommand{\A}{\mathscr A}
\newcommand{\B}{\mathscr B}
\newcommand{\x}{\times}
\newcommand{\<}{\langle}
\renewcommand{\>}{\rangle}
\renewcommand{\a}{\alpha}

\newcommand{\bs}{\bigskip}
\newcommand{\dt}{\cdot}
\newcommand{\C}{\mathbb C}
\newcommand{\ds}{\displaystyle}
\renewcommand{\d}{\delta}

\newcommand{\e}{\varepsilon}

\newcommand{\G}{\Gamma}
\renewcommand{\i}{\infty}
\renewcommand{\H}{\mathscr H}
\newcommand{\J}{\mathscr J}
\renewcommand{\l}{\lambda}

\newcommand{\M}{\mathscr M}
\newcommand{\mb}{\mbox}
\renewcommand{\ni}{\noindent}
\newcommand{\N}{\mathbb N}
\renewcommand{\O}{\Omega}
\renewcommand{\o}{\omega}
\newcommand{\ot}{\otimes}
\newcommand{\p}{\partial}
\newcommand{\R}{\mathscr R}
\newcommand{\Sig}{\Sigma}
\newcommand{\s}{\sigma}
\newcommand{\stt}{\subset}
\renewcommand{\t}{\tau}
\renewcommand{\th}{\theta}
\newcommand{\var}{\varphi}
\newcommand{{\z}}{\Bbb Z}
\newcommand{\pr}{probabili}
\newcommand{\ncc}{noncommutative }
\newcommand{\en}{entropy }
\newcommand{\dis}{distribution}
\newcommand{\mi}{microstates }
\newcommand{\jx}{{\mathscr J}(X:B)}

\begin{document}
\setlength{\baselineskip}{16pt}

\begin{center}
{\Huge Free Entropy}

\bs\bs {\sc Dan Voiculescu}

\bs Department of Mathematics\\University of California\\
Berkeley, California 94720--3840
\end{center}

\bs\bs\begin{center}
{\bf Abstract}
\end{center}

Free entropy is the analogue of entropy in free probability theory.  The paper is a
survey of free entropy, its applications to von Neumann algebras, connections to random
matrix theory and a discussion of open problems.

\bs\bs
\setcounter{section}{-1}\section{Introduction}

Entropy, from its initial appearance in thermodynamics and passing through
statistical mechanics and the mathematical theory of communications of Claude
Shannon, has come to play, in various guises, a fundamental role in many
parts of mathematics. This article is about a recent addition ([35]) to the
mathematical territory of \en.

Free \en refers to the analogue of \en in free \pr ty theory, i.e.  a
quantity playing the role of \en in a highly \ncc \pr stic framework, with
independence modelled on free products instead of tensor products. Free \pr
ty theory can be viewed as a parallel to some basic \pr ty theory drawn
starting from a new type of independence. Surprisingly, the parallelism of
the classical and free theories appears to go quite far, as illustrated for
instance by the existence of a free \en theory.

From another perspective, free \pr ty, and in particular free entropy, has
deep connections on one hand with the asymptotic behavior of large random
matrices and on the other hand with operator algebras. One consequence is
that the von Neumann algebras of free groups, once viewed as exotic
creatures, are now much better understood and perceived as important objects.

\bs
{\bf Acknowledgments.} Part of this work was done by the author
for the Clay Mathematics Institute. Partial support was also
 provided by
National Science Foundation grant DMS--0079945. 

\newpage
\section{Free Probability Background}

\subsection{Some basic laws}

Free \pr ty theory being a parallel to classical \pr ty theory, we may
compare the two by taking a look at corresponding fundamental distributions.
\begin{quote}
a) The role of the Gaussian distribution in free \pr ty theory is held by the
semi-circle distribution, which is a distribution with compact support.

b) For the Poisson distribution, the free correspondent is a distribution
related to the semi-circle law. It is also a compactly supported distribution
which has at most one atom.

c) The free Cauchy distribution is the Cauchy distribution itself, i.e. the
free correspondent is the same as the classical law.
\end{quote}

The semi-circle distribution occurs in random matrix theory, where Wigner
discovered that it is the limit distribution of eigenvalues of large
hermitian Gaussian matrices. Similarly, the free Poisson laws also occur in
random matrix theory as limit distributions of eigenvalues for matrices of
the form $X^*X$ where $X$ is a rectangular Gaussian matrix, it is the
Pastur-Marchenko distribution.

Like in many other situations, relations among \pr ty \dis s signal
structural connections, in this case a connection between free \pr ty theory
and random matrix theory.

\bs \ni {\bf Figure 1}\bs

Free \pr ty theory can be described as \ncc \pr ty theory endowed with the
definition of free independence. The next sections briefly explain the two
terms: \ncc \pr ty theory and free independence.

\subsection{Noncommutative \pr ty theory}

In classical \pr ty theory, numerical random variables are measurable
functions on a space of events $\O$ endowed with a \pr ty measure $\mu$, i.e.
a positive measure of mass one. The expectation $E(f)$ of a random variable
$f$ is the integral $\int f\,d\mu$.

Roughly speaking, \ncc \pr ty theory replaces the ring of numerical random
variables by a possibly \ncc algebra ${\A}$ over ${\C}$ with unit
$1\in{\A}$, which is endowed with a linear expectation functional
$\var:{\A}\to{\C}$, such that $\var(1)=1$. $({\A},\var)$ is a {\it
\ncc \pr ty space} and elements $a\in{\A}$ are {\it \ncc random variables.}

Often $({\A},\var)$ is an algebra of bounded operators on a Hilbert space
${\H}$ and the functional $\var$ is defined by a unit-vector
$\xi\in{\H}$, i.e., $\var(a)=\<a\xi,\xi\>$. Typically quantum mechanical
quantities $a\in{\A}$ are described in this way and $\xi$ is the state-vector.

The {\it \dis \ of a random variable} $a\in{\A}$ is the linear map
$\mu_a:{\C}[X]\to{\C}$ so that $\mu_a(P)=\var(P(a))$. The information encoded
in $\mu_a$ is the same as giving the collection of moments
$(\var(a^n))_{n>0}$. 

Similarly for a family $\a =(a_i)_{i\in I}$ of random variables in ${\A}$,
the \dis \ is the map $\mu_{\a}:{\C}\<X_i\mid i\in I\>\to{\C}$ so that
$\mu_{\a}(P)=\var(P((a_i)_{i\in I}))$, where ${\C}\<X_i\mid i\in I\>$ is the
algebra of \ncc polynomials in the indeterminates
$(X_i)_{i\in I}$. Like in the one-variable case, $\mu_\a$ contains the same
information as the \ncc moments $\var(a_{i_1}a_{i_2}\dots a_{i_p})$.

In the case of a self-adjoint operator $a=a^*$, \ $\mu_a$ can be identified
with a compactly supported \pr ty measure on $a$. Indeed if $E(a;\o)$ is the
projection-valued spectral measure of $a$, then
\[
\var(P(a)) =\<P(a)\xi,\xi\> =\int
P(t)\<dE(a;(-\i,t)\xi,\xi\>
\]
i.e., $\mu_a$ ``is" $E(a; \,\dt\,)\xi,\xi\>$.

The usual context for free entropy theory is the more restricted one of a 
{\it tracial $W^*$-\pr ty space} $(M,\t)$. This means that $M$ is a 
$W^*$-algebra (synonymous to von Neumann algebra) and that the expectation
function $\t$ is a trace. This means that $M$ is a  self-adjoint \alg \ of
bounded operators on a Hilbert space ${\H}$ (i.e. $T\in M\Rightarrow T^*\in
M$) which is weakly closed (i.e., if for some net $(T_i)_{i\in I}$ in $M$,
we have $\<T_ih,k\>\to \<Th,k\>$ for all $h,k\in{\H}$, then $T\in M$).
The condition on $\t$ is that $\t(ST)=\t(TS)$ for all $S,T\in M$.

If $(\O,\Sig,\mu)$ is a \pr ty space then $M=L^{\i}(\O,\Sig,\mu)$ acting as
multiplication operators on $L^2(\O,\Sig,\mu)$ is a $W^*$-\alg \ and the
expectation functional $\t$ defined by the vector $1\in L^2$ is trivially a
trace since $M$ is commutative. Note that $\t$ coincides with the classical
expectation functional on $L^{\i}$ defined by $\mu$. Thus tracial
$W^*$-\pr ty spaces subsume the context of classical \pr ty spaces.

A fundamental class of tracial $W^*$-\pr ty spaces  is generated by discrete
groups $G$. Let $\l$ be the left regular representation of $G$ on
$\ell^2(G)$, i.e.~$\l(g)e_h=e_{gh}$ where $e_g$, $g\in G$, are the canonical
basis vectors in $\ell^2(G)$. Then {\it the von Neumann \alg} \ $L(G)$ is
defined as the weakly closed linear space of $\l(G)$. Roughly speaking,
$L(G)$ consists of those left convolution operators $\sum_{g\in G} c_g\,\l(g)$
which are bounded on $\ell^2$. The trace $\t$ is the 
von Neumann trace which is defined by the basis vector $e_e$ (or any other
$e_g)$. Note that $\t(\sum_g c_g\,\l(g))=c_e$ (here the next $e$ denotes the
neutral element in $G$).

\subsection{Free independence}

A family of sub\alg s $({\A}_i)_{i\in I}$, with $1\in{\A}_i$, in
$({\A},\var)$ is {\it freely independent} if
\[
\var(a_1\dots a_n) \ = \ 0
\]
whenever $\var(a_j)=0$, $1\leq j\leq n$ and $a_j\in{\A}_{i(j)}$
with $i(j)\neq i(j+1)$, $1\leq j\leq n-1$. A family of subsets 
$(\o_i)_{i\in I}$ in $(A,\var)$ is freely independent if the \alg s
${\A}_i$ generated by $\{1\}\cup\o_i$ are freely independent.

The above definition means that products of centered variables, such that
consecutive ones are in different \alg s, have expectation zero. Note that
this does not preclude that $i(j)=i(k)$ as long as $|j-k|\geq 2$.

In general free independence requires that variables be very far from
commuting. For instance, if $X,Y$ are freely independent and centered
$\var(X)=\var(Y)=0$, then the free independence condition requires that
$\var(XYXY)=0$ while commutation of $X$ and $Y$ would imply
\[
\var(XYXY)=\var(X^2Y^2)=\var(X^2)\var(Y^2)
\]
where the last equality is derived from free independence
\[
\var((X^2-\var(X^2)1)(Y^2-\var(Y^2)1)) = 0 \ .
\]
Thus commutation is impossible if $\var(X^2)\neq 0$, \ 
$\var(Y^2)\neq 0$.

A basic example of free independence is provided by groups.
A family $(G_i)_{i\in I}$ of subgroups of a group $G$ is free, in the sense
of group theory if there is no non-trivial \alg ic relation in $G$ among the
$G_i$'s which translates into the requirement that $g_1g_2\dots g_n\neq e$
whenever $g_j\neq e$, $1\leq j\leq n$ and $g_j\in G_{i(j)}$ with
$i(j)\neq i(j+1)$, $1\leq j\leq n\!-\!1$. It can be shown that in 
$(L(G),\t)$ the free independence of the sets $(\l(G_i))_{i\in I}$ is
equivalent to the requirement that the family of subgroups 
$(G_i)_{i\in I}$ is \alg ically free. Note that this is also equivalent to
the \fr of the von Neumann \alg s generated by the $\l(G_i)$.

\subsection{Random matrices in the large $N$ limit}

The explanation found in [33] for the clues to a connection between free \pr
ty and random matrices is that \fr occurs asymptotically among large random
matrices.

Very roughly the connection is as follows. A \rn matrix is a classical
matrix-valued \rn variable. At the same time \rn matrices give rise to
operators, i.e.~to \ncc \rn variables. Note that the passage from the
classical variable to the \ncc one means forgetting part of the information
(the \ncc moments can be computed from the classical \dis \ but not vice
versa). Then under certain conditions (like unitary invariance) independent
\rn matrices give rise asymptotically as their size increases to freely
independent \ncc \rn variables.

The \ncc \pr ty framework for \rn matrices is given by the \alg s
\[
{\A}_N \ = \ L^{-\i}(X,{\mathscr M}_N)
\]
where $(X,\Sig,d\sigma)$ is a \pr ty space, ${\mathscr M}_N$ denotes the $N\x
N$ complex matrices and $L^{-\i}$ stands for the intersection  of
$L^p$-spaces $1\leq p <\i$. The expectation functional on ${\A}_N$ is 
$\var_N:{\A}_N\to{\C}$ given by
\[
\var_N(T) \ = \ N^{-1}\int_X {\mb{Tr}}(T(x))d\sigma(x) \ .
\]

The simplest instance of asymptotic \fr is provided by a pair of Gaussian
matrices. Let $T_j^{(N)} =(a^{(N)}_{p,q;j})_{1\leq p,q\leq N}\in{\A}_N$,
$j=1,2$, where $a^{(N)}_{p,q;j} =a^{(N)}_{q,p;j}$ and 
$\{a^{(N)}_{p,q;j}\mid 1\leq p\leq q\leq N, \ j\!=\!1,2\}$ 
are independent $(0,N^{-1})$-Gaussian. Then $T_1^{(N)},T_2^{(N)}$
are asymptotically free as $N\to \i$, in the sense that the \alg ic relations
among the \ncc moments of the pair $(T_1^{(N)},T_2^{(N)})$ which represent
the \fr conditions, are satisfied in the limit $N\to\i$.

Among the uses of asymptotic freeness of \rn matrices are the study of the
large $N$-limit of \rn matrices with free \pr ty techniques on one hand and
on the other hand the operator \alg \ applications. Operator \alg s such as
the von Neumann \alg s of free groups $L(F(n))$ are generated by free \rn
variables and can therefore be viewed as asymptotically generated by \rn
matrices. This has provided the intuitive background for many new results.

\subsection{Free independence with amalgamation}

In usual \pr ty theory conditional independence amounts to replacing the
scalar expectation functional with the conditional expectation w.r.t.~a
sub-$\s$-\alg \ of events, i.e., the expectation takes values in a sub-\alg \
of the \alg \ of \rn variables.

The free analogue of conditional independence is \fr with amalgamation.
The context is a ${\B}$-valued \pr ty space, i.e. $({\A},E,{\B})$
where $1\in{\B}\stt{\A}$ is an inclusion of unital \alg s over ${\C}$
and $E:{\A}\to{\B}$ is ${\B}\!-\!{\B}$-bilinear and
$E|_{\B}={\mb{id}}_{\B}$. Then a family of sub\alg s $({\A}_i)_{i\in I}$,
${\B}\stt{\A}_i\stt{\A}$ is ${\B}$-freely independent if $E(a_1\dots a_n)=0$
whenever $E(a_j)=0$, $1\leq j\leq n$, $a_j\in{\A}_{i(j)}$,
$i(k)\neq i(k+1)$, $i\leq k<n$.

If $(M,\t)$ is a tracial $W^*$-\pr ty space, with faithful $\t$
(i.e., $\t(x^*x)=0\Rightarrow x=0$) then there are canonical conditional
expectations onto von Neumann sub\alg s. If $I\in N\stt M$ is a von Neumann
sub\alg , then $\<m_1,m_2\>=\t(m^*_2m_1)$ is an inner product on $M$ and
$E_N$ is defined as the orthogonal projection of $M$ onto the Hilbert space
completion of $N$. It turns out that actually $E_N(M)\stt N$ and
$\|E_Nm\|\leq \|m\|$. Of course for the $L^2$-norm
$|m|_2=(\t(m^*m))^{\frac 12}$ we also have $\|E_Nm\|_2\leq \|m\|_2$.
Moreover $E_N$ is $N\!-\!N$-bilinear. This is clearly a generalization of the
classical situation where $M=L^{\i}(\O,\Sig,\mu)$ and
$N=L^{\i}(\O,\Sig_1,\mu)$ with $\Sig_1\stt\Sig$ a $\s$-sub\alg .

If $G$ is a group and $H$ a subgroup let $L(H)$ be identified with the
$W^*$-sub\alg \ generalized by $\l(H)$ in $L(G)$. Then
\[
E_{L(H)} \sum_{g\in G} c_g\, \l(g) \  = \
\sum_{g\in H} c_g\, \l(g) \ .
\]
Also if $H\stt G_i\stt G$ is a family of subgroups  indexed by $I$, then the
$L(G_i)$ are $L(H)$-freely independent in $(L(G),E_{L(H)})$ iff the subgroups
$G_i$ are \alg ically free with amalgamation over $H$.

\subsection{Background references}

The beginning of free \pr ty theory is the paper [31] and the
connection to \rn matrices is in [33]. A comprehensive introduction to 
free \pr ty theory is given in [42] and for \pr sts (i.e.~for readers who
prefer operator \alg s kept to a minimum) there are the St-Flour lectures
[39]. Some standard operator \alg \ books are [7], [8], [19], [29].

\section{Matricial Microstates Approach to Free Entropy}

\subsection{Underlying idea}

Shannon's \en of a continuous $n$-dimensional \dis \ ([23])
is given by the formula
\[
H(f_1,\dots ,f_n) =-\int_{{\R}^n} p(t_1,\dots ,t_n)\log p
(t_1,\dots ,t_n) dt_1\dots dt_n
\]
where $f_1,\dots ,f_n$ are real-valued \rn variables with Lebesgue absolutely
continuous joint \dis \ with density $p(t_1,\dots ,t_n)$.
A free analogue to $H(f_1,\dots ,f_n)$ will be a number
$\chi(X_1,\dots ,X_n)$ [35 II] associated to an $n$-tuple of self-adjoint
elements $X_j$ $(1\leq j\leq n)$ in a tracial $W^*$-\pr ty space
$(M,\t)$, the properties of $\chi$ w.r.t.~\fr being parallel to those of $H$
w.r.t.~classical independence. 

Information-theoretic and physical \en though different concepts, have also
much in common. In particular, the formula for $H(f_1,\dots ,f_n)$ can be
derived from Boltzmann's fundamental formula $S=k\log W$. The connection to
the Boltzmann formula and the fact that \fr occurs asymptotically among large
matrices, are the key to the definition of $\chi$.

Boltzmann's  formula says that the \en $S$ of a ``macrostate" is
proportional to the logarithm of its ``Wahrscheinlichkeit" $W$ (\pr ty),
where the \pr ty of the ``macrostate" is obtained by counting how many 
``microstates" correspond to that ``macrostate". For mathematical purposes,
microstates are often associated with a given degree of approximation,
and one then takes a normalized limit when the number of microstates goes to
infinity, followed by a limit improving the approximation.

For simplicity, here is how this works for the \en of a discrete \rn variable
with outcomes $\{1,\dots ,n\}$ with \pr ties $p_1,\dots ,p_n$.
The microstates are the set $\{1,\dots ,n\}^N=\{f\mid f:\{1,\dots ,N\}\to
\{1,\dots n\}\}$ and the microstates which approximate the discrete \dis \
are $\G(p_1,\dots ,p_n;\e,N)$ consisting of those $f$ such that
\[
\left| \frac{|f^{-1}(j)|}{N}-p_j\right| < \e
\]
($|f^{-1}(j)|$ the number of elements in the pre-image.) One then takes the
limit of
\[
N^{-1}\log|\G(p_1,\dots ,p_n;\e,N)|
\]
as $N\to\infty$ and then lets $\e$ go to zero. Using repeatedly Stirling's
formula one gets the familiar $-\Sig p_j\log p_j$ result in the end.

To define $\chi$, the microstates will be matricial.

\subsection{The definition of $\chi(X_1,\dots ,X_n)$ \ {\footnotesize{\rm [35
II]}}}

Given $X_j=X^*_j\in M$, $1\leq j\leq n$, where $(M,\t)$
is a tracial $W^*$-\pr ty space, the set of approximating matricial 
\mi will be denoted $\G_R(X_1,\dots , X_n;m,k,\e)$ where $R>0$,
$m\in{\N}$, \ $k\in {\N}$, $\e>0$. Here $R$ is a cut-off parameter, $k$ the
size of matrices and $(m,\e)$ the degree of approximation. With 
${\mathscr M}^{sa}_k$ denoting the self-adjoint complex $k\x k$ matrices,
the approximating \mi are $n$-tuples $(A_1,\dots ,A_n)\in ({\M}^{sa}_k)^n$
such that
\[
|\t(X_{i_1}\dots X_{i_p})-k^{-1}{\mb{Tr}}(A_{i_1}\dots A_{i_p})|<\e
\]
for all $1\leq p\leq m$, $(i_1,\dots ,i_p)\in \{1,\dots ,n\}^p$
and $\|A_j\| <R$, $1\leq j\leq m$.

Let \ vol \ denote the euclidean volume on $(M^{sa}_k)^n$ w.r.t.~the
Hilbert-Schmidt scalar product
\[
\<(A_1,\dots ,A_n),(B_1,\dots ,B_n)\>=\sum_j {\mb{Tr}} \ A_jB_j \ .
\]
Taking 
\[
\limsup_{k\to\i} (k^{-2}\log {\mb{vol}}\G_R (X_1,\dots ,X_n;m,k,\e)+
\frac n2\log k)
\]
and then
\[
\sup_{R>0} \ \inf_{m\in{\N}} \ \inf_{\e>0}
\]
of the result, we obtain \ $\chi(X_1,\dots ,X_n)$.

Note that the cut-off $R$ has only a minor influence, instead of the sup over
$R$ we could have taken a fixed $R$ larger than $\|X_j\|$,
$1\leq j\leq n$.

\subsection{Basic properties of $\chi(X_1,\dots ,X_n)$}

\begin{description}
\item[$\chi.1.$] {\bf Upper Bound} ([35 II])\\
$\chi(X_1,\dots ,X_n)\leq 2^{-1}n\log (2\pi en^{-1}C^2)$
where $C^2=\t(X^2_1+\dots +X^2_n)$.\\In particular
$\chi(X_1,\dots ,X_n)$ is either finite or $-\i$.

\item[$\chi.2.$] {\bf Subadditivity} ([35 II])\\
$\chi(X_1,\dots ,X_{m+n})\leq \chi(X_1,\dots ,X_m)+\chi
(X_{m+1},\dots ,X_{m+n})$.

\item[$\chi.3.$] {\bf Semicontinuity} ([35 II])\\
Assume $\|X_j^{(p)}\|\leq C<\i$, $1\leq j\leq n$, $p\in{\N}$
and $(X_1^{(p)},\dots ,X_n^{(p)})$\\ converges in \dis \ to 
$(X_1,\dots ,X_n)$, i.e.
\[
\lim_{p\to \i} \t\big(X^{(p)}_{i_1}\dots X^{(p)}_{i_k}\big) =\t
(X_{i_1}\dots X_{i_k})
\]
for all \ncc moments. Then
\[
\limsup_{p\to \i} \chi(X^{(p)}_1,\dots ,X^{(p)}_n)\leq
\chi(X_1,\dots ,X_n) \ .
\]

\item[$\chi.4.$] {\bf One Variable Case.} ([35 II])\\
$
\chi(X)=\iint\log |s-t|d\mu(s)d\mu(t) +\textstyle{\frac 34+\frac 12}
\log 2\pi$\\ 
where $\mu$ denotes the \dis \ of $X$. Thus, up to constants, $\chi(X)$
is\\minus the logarithmic energy of $\mu$.

\item[$\chi.5.$] {\bf Additivity and Free Independence.}  ([35 IV])\\
Assume $\chi(X_j)>-\infty$, \ $1\leq j\leq n$. Then\\
$\chi(X_1,\dots ,X_n)=\chi(X_1)+\dots +\chi(X_n)$ {\bf iff}
$X_1,\dots ,X_n$ are freely independent.

\item[$\chi.6.$] {\bf Semicircular Maximum.} ([35 II])\\
Assume $\t(X^2_1)=\dots=\t(X^2_n)=1$. Then $\chi(X_1,\dots ,X_n)$ is maximum
{\bf iff} $X_1,\dots ,X_n$ are freely independent and have (0,1)-semicircular
\dis s.

\item[$\chi.7.$] {\bf Infinitesimal Change of Variables.} ([35 IV])\\
Let ${\C}\<t_1,\dots ,t_n\>$ be the ring of \ncc polynomials in the
indeterminates $t_1,\dots ,t_n$ endowed with the involution \ * \ so that
$(ct_{i_1}\dots t_{i_p})={\bar c}t_{i_p}\dots t_{i_1}$. Then:
\begin{eqnarray*}
&&\left.\frac{d}{d\e}\,\chi(X_1+\e P_1(X_1,\dots ,X_n),\dots ,X_n+
\e P_n(X_1,\dots ,X_n))\right|_{\e=0} \\ 
&&\qquad= \sum_{1\leq j\leq j}
(\t\otimes\t)(\p_j P_j(X_1,\dots ,X_n))
\end{eqnarray*}
where $P_j=P^*_j\in{\C}\<t_1,\dots ,t_n\>$ and 
$\p_j:{\C}\<t_1,\dots ,t_n\>\to {\C}\<t_1,\dots ,t_n\>\otimes
{\C}\<t_1,\dots ,t_n\>$ is given by $\p_j\dt t_{i_1}\dots t_{i_p}=
\sum_{i_k=j} t_{i_1}\dots t_{i_{k-1}}\otimes t_{i_{k+1}}\dots t_{i_p}$.

\item[$\chi.8.$] {\bf Degenerate Convexity.} ([35 III])\\
Assume $n\geq 2$ and there are trace-states $\t',\t''$ on
$A=W^*(X_1,\dots ,X_n)$ so that $\t'\neq \t''$ and $\t=\th\t'+(1-\th)t''$
on $A$, $0<\th <1$. Then \ $\chi((X_1,\dots ,X_n)=-\i$.
\end{description}

\bs\ni{\bf Remarks:}
\begin{description}
\item[a)] It is an important open problem, whether replacing the lim\,sup in
the definition of $\chi$ by a lim\,inf (as $k\to\i$) yields the same
quantity. While this is unresolved it is sometimes convenient to use
$\chi_{\o}$, $\o$ an ultrafilter on ${\N}$, the quantity obtained by
replacing the lim\,sup by a limit as $k\to\o$.
\item[b)] Generalizing the ``if\,"-part of $\chi.5.$~to groups of
variables runs into the problem discussed in a). There is a partial
generalization ([36])
\[
\chi_{\o}(X_1,\dots ,X_{m+n}) =\chi_{\o}(X_1,\dots ,X_m)+\chi_{\o}(X_{m+1}
,\dots ,X_{m+n})
\]
if $\{X_1,\dots ,X_m\}$ and $\{X_{m+1},\dots ,X_{m+n}\}$ are freely
independent.
\item[c)] We preferred to state the weaker infinitesimal version of the
change of variable formula because it is easier to state and will be used
later. Roughly the change of variable formula is of the form:
\[
\chi(F_1(X_1,\dots ,X_n),\dots , F_n(X_1,\dots ,X_n))=
\chi(X_1,\dots ,X_n)+\log |{\mb{det}}|(DF(X_1,\dots ,X_n))
\]
where there is a long  list of details about the \ncc power series
$(F_1,\dots ,F_n)$, the Kadison-Fuglede determinant $|$det$|$ and the
differential $DF$ for which the reader is referred to the original
paper [35 II].
 
\item[d)] Given $X_1,\dots ,X_n$ and $m\in{\N}$, $\e>0$ is there
$k\in{\N}$ and $R>0$ so that
\[
\G_R(X_1,\dots ,X_n)m,k,\e)\neq\emptyset \ ?
\]
This very basic question is equivalent to a problem of A.Connes on embedding
II${}_1$-factors into the ultraproduct of the hyperfinite II${}_1$-factor.

\item[e)] The ``if\,"-part of $\chi.5.$~relies essentially on
asymptotic freeness of \rn matrices. What the result and its proof show, is a
sharp difference between one- and multi-\rn matrix theory. Roughly, if $n=1$,
then sets of \mi $\G_R(X;m,k,\e)$ will be like tubes around
the unitary orbit of some microstate $\{UAU^*\mid U\in{\mathscr U}(n)\}$.
If $X_1,\dots ,X_n$ are freely independent and $n>1$, then
$\G_R(X_1,\dots ,X_n;m,k,\e)$ is much larger than a tube around
$\{(UAU^*,\dots ,UA_nU^*)\mid U\in{\mathscr U}(n)\}$, actually up to sets,
the measure of which goes to 0 as $k\to\i$, it is more like the product of
tubes around the orbits of the components, i.e. $\{UA_kU^*\mid
 U\in{\mathscr U}(n)\}$.
\end{description}

\subsection{The free entropy dimension\ {\footnotesize{\rm [35
II]}}}

The free \en being a normalized limit of logarithms of volumes of sets of
matricial microstates, there is also a corresponding normalized dimension of
sets of microstates. The definition is reminiscent of the definition of the
Minkowski content.

The {\it free \en dimension} $\d(X_1,\dots ,X_n)$ is given by the formula
\[
\d (X_1,\dots ,X_n) \ = \ n+\limsup_{\e\downarrow 0} \
\frac{\chi(X_1+\e S_1,\dots ,X_n+\e S_n)}{|\log \e|}
\]
where $S_1,\dots ,S_n$ have (0,1)-semicircular \dis s and
$\{X_1,\dots ,X_n\},\{S_1\},\dots ,\{S_n\}$ are freely independent.

In a  number of applications it is necessary for technical reasons 
to use a modification $\d_0(X_1,\dots ,X_n)$ of $\d$. It is not known whether
$\d$ and $\d_0$ are actually different. $\d_0$ is obtained by replacing
$\chi(X_1+\e S_1,\dots ,X_n+\e S_n)$ in the definition of $\d$ by
$\chi(X_1+\e S_1,\dots ,X_n+\e S_n:S_1,\dots ,S_n)$ where
$\chi(X_1,\dots ,X_n:Y_1,\dots ,Y_p)$ is defined like $\chi$ using
\[
\G_R(X_1,\dots ,X_n:Y_1,\dots ,Y_p;m,k,\e) =
pr_{\{1,\dots ,n\}}
\G_R(X_1,\dots ,X_n,Y_1,\dots ,Y_p;m,k,\e)
\]
Since all this becomes rather technical, we will limit our discussion to $\d$
in the rest of this section.

Here are some basic properties of $\d$.
\begin{quote}
{\bf a)} $\d(X_1,\dots ,X_n)\leq n$. We also have 
$\d(X_1,\dots ,X_n)\geq 0$ when the problem in\\
\mbox{}\qquad 2.3--Remark d) has an
affirmative answer for $X_1,\dots ,X_n$.\\
{\bf b)} $\d( X_1,\dots ,X_{m+n})\leq \d(X_1,\dots ,X_m)+
\d(X_{m+1},\dots ,X_{m+n})$\\
{\bf c)} $\d(X_1,\dots ,X_n)=\d(X_1)+\dots +\d(X_n)$ if $X_1,\dots ,X_n$
are freely independent.\\
{\bf d)} $\d(X)=1-\ds{\sum_{t\in{\R}}}(\mu(\{t\}))^2$ where $\mu$ is the
\dis \ of $X$.\\
{\bf e)} $\chi(X_1,\dots ,X_n) >-\i\Rightarrow \d(X_1,\dots ,X_n)=n$.
\end{quote}

\subsection{Operator algebra applications}

Free \en has led to new results on von Neumann algebras, in particular the
solution of some old problems has been found. The new results are about
separable II${}_1$ factors, i.e. von Neumann algebras $M$ of infinite
dimension acting on separable Hilbert spaces, which have a faithful
trace-state $\t$ and trivial center $Z(M)={\C}I$. Typical examples are the
$L(G)$'s where $G$ is a countable discrete group with infinite conjugacy
classes.

\begin{description}
\item[$1^{\circ}$] {\bf Absence of Cartan Subalgebras} ([35 III])\\
The free group factors $L(F(n))$ $(n\geq 2)$ have no  Cartan subalgebras.
A Cartan subalgebra $A\stt M$ ($M$ a II${}_1$ factor) is a maximal abelian
$W^*$-subalgebra, the normalizer of which $N(A)=\{u\in M\mid u$ unitary,
$uAu^*=A\}$ generates $M$. The concept mimics the properties of the algebra
of diagonal matrices inside the algebra of $n\x n$ matrices. $M$ has a 
Cartan subalgebra iff it  can be obtained from an ergodic measurable
equivalence relation via a construction of Feldman and Moore ([13]).
It was an open problem whether all separable II${}_1$ factors arise this way
from ergodic theory.

\item[$2^{\circ}$] {\bf Prime II${}_1$ factors} ([15 II])\\
$L(F(n))$ $(n\geq 2)$ is prime, i.e. is not a 
$W^*$-tensor product $M_1\otimes M_2$ of $\i$-dimensional von Neumann
algebras. 
The existence of separable II${}_1$ factors was also an old open question.

\item[$3^{\circ}$] {\bf Products of abelian subalgebras} ([30])\\
If $n$ is large enough, $L(F(n))$ is not the 2-norm closure of the linear
span of a product $A_1\dots A_m$ of $m$ abelian $W^*$-subalgebras.

Using a fundamental theorem of A.Connes, by which all separable
II${}_1$-factors $L(G)$ with $G$ amenable are isomorphic, it follows that
in the amenable case $L(G)=\overline{\mb{span}\, A_1A_2}$ for a pair of
  abelian $W^*$-subalgebras. This is in sharp contrast with the $L(F(n))$
situation.
\end{description}

The principle underlying the proofs of these results is to show that a
certain property (existence of a Cartan subalgebra, non-primeness, product of
abelian, etc.) implies that a generator $X_j=X^*_j$ $(\lejn)$ of the von
Neumann algebra has $\chi(X_1,\dots ,X_n)=-\i$. On the other hand
$L(F(n))$ has a generator $X_1,\dots ,X_n$ with $\chi(X_1,\dots ,X_n)>-\i$
(consider Borel-logarithms of the generating unitaries $\l(g_1)\dots\l(g_n)$
and use $\chi.4$ and  $\chi.5$).
This kind of result, started by the absence of Cartan algebras result ([35
III]) has meant developing increasingly ingenious ways of estimating
volumes of matricial \mi for generators ([12],[15],[30]). 

Note also that for most of the above results there are stronger forms,
where $\chi(X_1,\dots ,X_n)=-\i$ is replaced by $\d_0(X_1,\dots
,X_n)\leq 1$ for a generator. In this direction there is also the following
recent result.

\begin{description}
\item[$4^{\circ}$] {\bf Property $T$} ([15 III])\\
If $X_j=X^*_j$ $(\lejn)$ is a generator of $L(SL(rm+1;{\mathbb Z}))$ $(m\geq
1)$ then $\d_0(X_1,\dots ,X_n)\leq 1$.
\end{description}
The restriction to odd numbers $2m+1$ is only to insure
factoriality (i.e.~trivial center).

\subsection{Comments on the \mi approach}

The use of \mi, per se, in the definition of free entropy,
should not bother us too much. There are many other situations in
mathematics where huge auxiliary objects are used to define some
basic invariants (singular homology may come to mind for instance).
On the other hand,  the technical difficulties in this approach which
prevented us from completing the theory (see, for instance, Remarks a)
and b) in 2.3) are a problem.

Much impetus for further developing free \en theory is provided by
von Neumann algebras. There is some hope that with stronger free \en
tools at hand, the currently best known problem in the area may be
settled in the affirmative:

\begin{quote}{\bf $L(F(n))$ isomorphism problem.}\\
Does $L(F(n))\simeq L(F(m))$ \ imply \ $m=n$\,?\end{quote}

\ni An even more far-fetched question is whether for the free \en
dimension, or for some variant of it, there is an affirmative answer
to:

\begin{quote}{\bf The \en dimension problem.}
If \ $X_j=X^*_j\in M$, \ $Y_k=Y^*_k\in M$, \ $\lejn$, \  $1\leq k\leq
n$, \ does \ 
$W^*(X_1,\dots ,X_n)=W^*(Y_1,\dots ,Y_m)$ \ imply \
$\d (X_1,\dots ,X_n)=\d(Y_1,\dots ,Y_m)$\,?\end{quote}

\ni Under certain conditions, an affirmative answer to the preceding
problem would follow (see [35 II]) from an affirmative answer to:

\begin{quote}{\bf Semicontinuity of $\d$ problem.}
If \ $X_j^{(p)}=X_j^{(p)^*}\in M$, \ $X_j=X^*_j\in M$, $\lejn$, \
$p\in{\N}$ are so that $s-\ds{\lim_{p\to\infty}}X_j^{(p)}=X_j$
 \ does it follow that 
$\ds{\liminf_{p\to\infty}\,\d (X_1^{(p)},\dots
,X_n^{(p)})\geq\d(X_1,\dots ,X_n)}$\,?\end{quote}

Little is known about these questions. About the semicontinuity
problem it is only known that in the rather uninteresting case $n=1$,
the answer is yes ([35 II]). For certain variants of $\d$, the much
weaker free \en dimension problem, with the $W^*$-algebras replaced
by the \alg s (no closures) of the $X$'s and $Y$'s, the answer is
affirmative ([36]). Also the isomorphisms of various free product von
Neumann \alg s ([10],[11],[22],[34]) seem not to contradict the
invariance of $\d$ on generators. Finally, it is known [22] that
there are only two possibilities in the isomorphism problem: either
all $L(F(n+1))$, $n\in {\N}\cup\{\i\}$ are isomorphic or all are
non-isomorphic.

\section{Infinitesimal Approach to Free Entropy}

\subsection{Fisher information background}

The Fisher information ${\J}(f)$ of a real \rn variable $f$ is the
derivative of the \en in the direction of a Brownian motion starting
at $f$, or equivalently:
\[
\textstyle{\frac 12}{\J}(f) = \lim_{\e\downarrow 0}
(H(f+\e^{\frac 12}g)-H(f))
\]
where $g$ is a (0,1)-Gaussian variable independent of $f$.
Using the Brownian motion starting at $f$ one can then express $H$
via ${\J}$,
\[
H(g)-H(f) =\textstyle{\frac 12}\int^\i_0 ({\J}(f+t^{\frac 12}g)
-(1+t)^{-1})dt \ .
\]

On the other hand, if the \dis  \ of $f$ is Lebesgue absolutely
continuous with smooth density $p$, then one finds
\[
{\J}(f) =\int_{\R} \frac{(p'(t))^2}{p(t)} \ dt \ .
\] 
The last formula can also be expressed as an $L^2$-norm
\[
{\J}(f) = \left\| \frac{p'}{p}\right\|^2_{L^2({\R},pd\l)}
\]
or equivalently
\[
{\J}(f) = E\left(\big(\frac{p'}{p}(f)\big)^2\right) \ .
\]
The Fisher information initially appeared in statistics, where it was
defined by the preceding formula with $\frac{p'}{p}\, (f)$
being the so-called score-function of $f$. The score is also
fundamental for other reasons:
a) infinitesimally the effect on the \dis s of the perturbations
$f+\e^{\frac 12}g$ and $f+\frac \e 2 \ \frac{p'}{p}\,(f)$ is the same;
b) the score is a gradient for the \en when the perturbations of $f$
are of the form $f+\e Q(f)$ where $Q$ is a polynomial.

Related to property a) of the score the element
$p'/p\in L^2({\R},pd\l)$ can also be described as:
\[
 \frac{p'}{p} \ = \ -\big( \frac{d}{dt}\big)^*\,1
\]
where $\frac{d}{dt}$ is the operator of derivation densely defined on
polynomials in $L^2({\R},pd\l)$ and $p$ is smooth with compact
support. In particular,
\[
{\J}(f) = \left\|
\big(\frac{d}{dt}\big)^*\,1\right\|^2_{L^2({\R},pd\l)}
\]

Based on properties of the free \en $\chi$ and on one-dimensional
computations [35 I], it turns out [35 V] that the free analogue of
the Fisher information can be obtained, roughly speaking, by
replacing the operator of derivation $d/dt$ by some difference
quotient, which sends a polynomial $P(t)$ to the two-variable
polynomial:
\[
\frac{P(s)-P(t)}{s-t} \ .
\]
Dealing with several  noncommuting variables will involve \ncc
generalizations of the difference quotient, like the derivations
appearing in the infinitesimal change of variable formula for $\chi$.

\subsection{The free difference quotient}

Let $X=X^*\in M$ and $1\in B\stt M$ be a $*$-sub\alg \ such that $X$
and $B$ are \alg ically free (i.e., no non-trivial \alg ic relation
between $B$ and $X$). We denote by $B[X]$ the \alg \ generated by $B$
and $X$ and consider the linear map:
\[
\p_{X:B}: B[X] \to B[X]\otimes B[X]
\]
so that 
\[
\p_{X:B} b_0Xb_1 X\dots b_n \ = \
\sum_{1\leq k\leq n} b_0X\dots b_{k-1}\otimes b_kX\dots b_n \ .
\]
With the natural $B[X]$-bimodule structure on $B[X]\otimes B[X]$,
the map $\p_{X:B}$ is a derivation and it is the only one such that
$\p_{X:B}|B=0$ and $\p_{X:B}X=1\otimes 1$.

Note that the partial derivation appearing in the infinitesimal
change of variable formula for $\chi$ correspond to taking
$B={\C}[X_1,\dots ,{\widehat X_j},\dots ,X_n]$ and
$X=X_j$ (here $X_1,\dots ,X_n$ are \alg ically free, noncommuting).

$B[X]$ is a linear subspace of $L^2(M,\t)$ and we shall consider
$L^2(B[X],\t)$ the closure of $B[X]$.

\subsection{The conjugate variable ${\J}(X:B)$ \ {\footnotesize{\rm [35 V]}}}

In the context of the preceding section $\p_{X:B}$ is a densely defined
unbounded operator from $L^2(B[X],\t)$ to $L^2(B[X],\t)\ot L^2(B[X],\t)$.
We define ${\J}(X:B)=\p^*_{X:B}1\ot 1$ if it exists and call it the conjugate
variable to $X$ (w.r.t.~$B$).

Several other names are appropriate for ${\J}(X:B)$: \ncc Hilbert transform,
free Brownian gradient, free score. All these designations correspond to
properties of ${\J}(X:B)$ which will be described in what follows. In
particular the passage from the usual (partial) derivative to the free
difference quotient justifies the ``free score" name.

Here are some basic facts about ${\jx}$.
\begin{description}
\item{${\mathscr J}\!.1.$ {\bf Hilbert transform.}}
If the distribution of $X$ is Lebesgue absolutely continuous and has
density $p\in L^3({\R},d\l)$, then ${\J}(X:{\C})=g(X)$, where
$g=2\pi Hp$, with $H$ denoting the Hilbert transform.

\item{${\mathscr J}\!.2.$ {\bf Enlarging the scalars.}}
If \ $1\in C\stt M$ is a $*$-sub\alg \ and $C$ and $B[X]$ are freely
independent in $(M,\t)$ then
\[
\jx = {\J}(X:B \vee C)
\]
where $B\vee C$ is the \alg \ generated by $B$ and $C$.
(There is a strengthening of this in [25]: it suffices to assume $C$ and
$B[X]$ are freely independent over $B$ in $(M,E_B)$.)

\item{${\mathscr J}\!.3.$ {\bf Semicircular perturbations.}}
If $S$ is (0,1) semicircular and $B[X]$ and $S$ are freely independent and
$\e >0$, then
\[
{\J}(X+\e S:B)b = \e^{-1} \, E_{B[X+\e S]}S \ .
\]
In particular, $\|{\J}(X+\e S:B)\|\leq 2\e^{-1}$, and the set of selfadjoint
$X$ for which  $\|{\J}(X:B)\| <\i$ is norm-dense in the selfadjoint part of
$M$.

\item{${\mathscr J}\!.4.$ {\bf Closability.}} If \ $|{\jx}|_2 <\i$ then
$\p^*_{X:B}$ is densely defined and $\p_{X:B}$ is closable.

\item{${\mathscr J}\!.5.$ {\bf Free Brownian gradient.}}
If $S$ is (0,1) semicircular, $B[X]$ and $S$ freely independent,
$|{\jx}|_2 <\i$ and $\e >0$, then:
\begin{eqnarray*}
&&\t(b_0(X+\frac\e 2 {\jx})b_1 (X+\frac\e 2 {\jx})\dots b_n) \\ [0.125in]
&&\qquad =\t(b_0(X+\e^{\frac 12}S)b_1(X+\e^{\frac 12}S)\dots b_n)+ O(\e^2) \ .
\end{eqnarray*}

\item{${\mathscr J}\!.6.$ {\bf Gradient of $\chi$.}}  Let 
$X_j=X^*_j\in M$, $\lejn$ and assume that\\$\chi(X_1,\dots ,X_n)>-\i$ and that 
${\J}_k={\J}(X_k:{\C}[X_1,\dots {\widehat X_k},\dots ,X_n])$, $1\leq k\leq n$
exist. Then
\[
\frac{d}{d\e}\, \chi(X_1+\e P_1,\dots,X_n+\e P_n)|_{\e=0} \ = \
\sum_{1\leq k\leq n} \t(P_k{\J}_k)
\]
where $P_k=P^*_k\in {\C}[X_1,\dots ,X_n]$, $1\leq k\leq n$.
\end{description}

\subsection{$\Phi^*(X_1,\dots ,X_n:B)$ \ {\footnotesize{\rm ([35 V])}}}
In the infinitesimal approach, the relative Fisher information
$\Phi^*(X_1,\dots ,X_n:B)$ of an $n$-tuple of selfadjoint variables 
$X_1,\dots ,X_n$ with respect to the sub\alg \ $B$ is defined by
\[
\Phi^*(X_1,\dots ,X_n:B) = \sum_{1\leq k\leq n}
|{\J}(X_k:B [X_1,\dots ,{\widehat X_k},\dots ,X_n])|^2_2
\]
if the right-hand side is defined and $+\i$ otherwise. The asterisk is to
distinguish quantities in this approach from the corresponding quantities in
the matricial \mi approach. 

Here are some properties of  $\Phi^*$.
\begin{description}
\item{$\Phi^*.1.$ {\bf Superadditivity.}}\\
$\Phi^*(X_1,\dots ,X_n,Y_1,\dots ,Y_m:B)\geq
\Phi^*(X_1,\dots ,X_n:B)+\Phi^*(Y_1,\dots ,Y_m:B)$

\item{$\Phi^*.2.$ {\bf Free additivity.}} If \
$B[X_1,\dots ,X_n]$ and $C[Y_1,\dots ,Y_m]$ are freely independent, then\\
$\Phi^*(X_1,\dots ,X_n,Y_1,\dots ,Y_m:B\vee C)=
\Phi^*(X_1,\dots ,X_n:B)+\Phi^*(Y_1,\dots ,Y_m:C)$.

\item{$\Phi^*.3.$ {\bf Free Cramer-Rao inequality.}}
$\Phi^*(X_1,\dots ,X_n:B)\t(X^2_1+\dots +X^2_n)\geq n^2$.
Equality holds iff $X_j$ are semicircular with $\t(X_j)=0$
$(\lejn)$ and $B,\{X_1\},\dots ,\{X_n\}$ are freely independent.

\item{$\Phi^*.4.$ {\bf Free Stam inequality.}}
If $B[X_1,\dots ,X_n]$ and $C[Y_1,\dots ,Y_m]$ are 
freely independent, then
\[
(\Phi^*(X_1+Y_1,\dots ,X_n+Y_n:B\vee C))^{-1}\geq
(\Phi^*(X_1,\dots ,X_n:B))^{-1} + (\Phi^*(Y_1,\dots ,Y_m:C))^{-1} \ .
\]

\item{$\Phi^*.5.$ {\bf Semicontinuity.}}
If \ $X_j^{(k)}=X_j^{(k)^*}\in M$ \ and \
$s-\ds{\lim_{k\to\i}} X_j^{(k)}=X_j$, then\\
$\ds{\liminf_{k\to\i}}\, \Phi^*(X_1^{(k)},\dots ,X_n^{(k)}:B)\geq
\Phi^*(X_1,\dots ,X_n:B)$.

\item{$\Phi^*.6.$}
If \ $\Phi^*(X_1,\dots ,X_n:B)=\Phi^*(X_1,\dots ,X_n:{\C})<\i$ then
$\{X_1,\dots ,X_n\}$ and $B$ are freely independent. If
$\Phi^*(X_1,\dots ,X_n,Y_1,\dots ,Y_m:{\C})=
\Phi^*(X_1,\dots ,X_n:{\C})+\Phi^*(Y_1,\dots ,Y_m:{\C})$ then
$\{X_1,\dots ,X_n\}$ and $\{Y_1,\dots ,Y_m\}$ are freely independent.
\end{description}

\subsection{$\chi^*(X_1,\dots ,X_n:B)$}
The free \en of $X_1,\dots ,X_n$ relative $B$, in the infinitesimal approach
is defined by
\[
\chi^*(X_1,\dots ,X_n:B) =
\textstyle{\frac 12}\int^\i_0
\left(\frac{n}{1+t} -\Phi^*(X_1+t^{\frac 12}S_1,\dots ,X_n+
t^{\frac 12}S_n:B\right)dt +\frac n2\log 2\pi e
\]
where the $S_j$'s are (0,1)-semicircular and
$B[X_1,\dots ,X_n],\{S_1\},\dots ,\{S_n\}$ are freely independent.

Here are some properties of  $\chi^*$.
\begin{description}
\item{$\chi^*.1.$} \ $\chi(X:{\C}) = \chi(X)$.
\item{$\chi^*.2.$} \ $\chi^*(X_1,\dots ,X_n)\leq\frac n2\log (2\pi n^{-1}C^2)$
\ where \ $C^2=\t(X^2_1+\dots +X^2_n)$.
 \item{$\chi^*.3.$} \ If $B[\xr]$ and $C$ are freely independent, then\\
$\chi^*(\xr:B)=\chi^*(\xr:B\vee C)$.
 \item{$\chi^*.4.$} {\bf Subadditivity.}\\
$\chi^*(\xr,Y_1,\dots ,Y_m:B\vee C)\leq
\chi^*(\xr:B)+\chi^*(Y_1,\dots ,Y_m:C)$
 \item{$\chi^*.5.$} {\bf Free additivity.} If \
$B[\xr]$ and $C[Y_1,\dots ,Y_m]$ are freely independent then the
inequality $\chi^*.4$ is an equality.
\item{$\chi^*.6.$} {\bf Semicontinuity.} If \ $\ds{s-\lim_{k\to\i}}
X^{(k)}_j=X_j$ \ then \\
$\ds{\limsup_{k\to\i}} \chi^*(X_1^{(k)},\dots ,X_n^{(k)}:B)\leq
\chi^*(\xr:B)$.
\item{$\chi^*.7.$} {\bf Information log-Sobolev inequality.}
If \ $\Phi^*(\xr:B)<\i$ then\\
$\chi^*(\xr:B)\geq \frac n2\log\ds{\left(
\frac{2\pi n e}{\Phi^*(\xr:B)}\right)}$, in particular\\
$\chi^*(\xr:B)>-\i$.
\end{description}

\subsection{Mutual free information and the derivation $\d_{A:B}$
\ {\footnotesize{\rm [35 VI]}}}
In the classical context, if $f,g$ is a pair of numerical \rn variables with
$H(f),H(g),H(f,g)$ finite, then their {\it mutual information} is
\[
I(f;g) \ = \ H(f) + H(g)-H(f,g) \ .
\]
Via an approximation procedure, the definition of $I(f,g)$ can be extended
well beyond the case of finite entropies (even Lebesgue absolute continuity of
distributions is not a requirement, see [6]). It also turns out that 
 $I(f,g)$  depends  only on the position of the von Neumann \alg s of $f$ and
$g$ inside the von Neumann algebra of $\{f,g\}$ endowed with the expectation
functional [in classical terms: the  triple of $\sigma$-\alg s of $f$,
respectively $g$, and respectively $(f,g)$-measureable events and the \pr ty
measure]. Note however that there is no infinitesimal theory for
$I(f,g)$ unless one is in the finite \en case and uses the infinitesimal
theory for entropy, i.e., there is no  infinitesimal
theory at the level of \alg s, since there is no natural deformation of the
pair of \alg s in sight.

In the free context, the situation is different. Given two von Neumann sub\alg
s $1\in A$, $1\in B$ in $(M,\t)$ there is a natural ``liberation process"
which deforms the pair $(A,B)$ to a freely independent pair:
$A,U(t)BU(t)^*$ where $\{U(t)\}_{t\geq 0}$ is a multiplicative unitary free
Brownian motion which is freely independent from $A\vee B$. 
This means $\{U(t)\}_{t\geq 0}$ is the free analogue of the corresponding
classical Brownian motion on the unit circle and can also be described, in
view of the asymptotic freeness of \rn matrices as the large $N$ limit of 
Brownian motions on the unitary groups $U(N)$ (see [2]). Via some heuristic
considerations this leads to an infinitesimal approach to a quantity
$i^*(A,B)$ which should play the role of the mutual free information for the
pair $(A,B)$.

The  infinitesimal approach  relies on a derivation
\[
\d_{A:B}:A\vee B \to (A\vee B)\ot (A\vee B)
\]
which exists under the assumption that $A$ and $B$ are \alg ically free (i.e.,
no non-trivial \alg ic relation). Here $(A\vee B)\ot (A\vee B)$ is with the
obvious $A\vee B$ bimodule structure and
\begin{center}
\begin{tabular}{lll}
$\d_{A:B} a$ &= \  $a\ot 1-1\ot a$ & if \ \ $a\in A$\\
$\d_{A:B} b$ &= \  0 & if \ \ $b\in B$ .
\end{tabular}
\end{center}

Like in the infinitesimal approach to free entropy, the key construction is
the liberation gradient 
\[
j(A:B) =\d^*_{A:B}1\ot 1
\]
where $\d_{A:B}$ is viewed as an unbounded operator densely defined
on $L^2(W^*(A\vee B),\t)$ with values in $L^2(W^*((A\vee B)\ot (A\vee
B)),\t\ot\t)$.

We list some of the main properties of $j(A:B)$.
\begin{description}
\item[j.1.] {\bf Liberation gradient.} \
$j(A:B) = -j(A:B)^*$ \ and \\
$\t\ds{\left(\prod^\to_{1\leq k\leq n}\! U(\e)a_kU(\e)^*b_k\right)} =
\t\ds{\left(\prod^\to_{1\leq k\leq n} \exp
\big(\frac{\e}{2} j(A:B)\big)a_k\,\exp\big(-\frac{\e}{2}
j(A:B)\big)b_k\right)} + O(\e^2)$\medskip\\
where $a_k\in A$, \ $b_k\in B$, $\ds{\prod^\to}$ denotes 
the ordered product and $(U(t))_{t\geq 0}$ is the multiplicative unitary free
Brownian motion free w.r.t.~$A\vee B$.
\item[j.2.] \ $j(A:{\C})=0$ \ and \
  $\ds{\sum_{1\leq k\leq n}} j(A_k:A_1\vee\dots\vee
A_{k-1}\vee A_{k+1}\vee\dots A_n)=0$
\item[j.3.] \ If $A,B,C$ is freely Markovian (i.e.~$A$ and $C$ are 
freely independent over $B$ in $(M,E_B)$) then
\begin{eqnarray*}
j(A:B) &=& j(A:B\vee C)\\
j(A:C) &=& E_{A\vee C}j(B:C)
\end{eqnarray*}
\item[j.4.] \ If $U$ is unitary and $A\vee B$ and
$\{U,U^*\}$ are freely independent, then
\[
j(A:UBU^*) = E_{A\vee UBU^*}\, j(A:B)
\]
and if the distribution of $U$ is absolutely continuous w.r.t.~Haar measure, 
$d\mu =pd\theta$, $p\in L^3$, then
\[
j(A:UBU^*) = -iE_{A\vee UBU^*}\, g(U)
\]
where \ $g(e^{i\th_1}) =\ds{-\,\frac{1}{2\pi} \ {\mb{p.v.}}\!\int
\frac{p(e^{i(\th_1-\th)})}{\tan(\th/2)}}\,d\th$ \
is the Hilbert transform.
\item[j.5.] \ 
$\ds{\|(E_A-E_{{\C}1})(E_B-E_{{\C}1})\|\leq
\frac{\|j(A:B)\|}{(1+\|j(A:B)\|^2)^{\frac 12}}}$\\
(The left-hand side is the norm of an operator on $L^2(M,\t)$.)
\item[j.6.] \
$j({\C}[\xr ]:B)=\sum_k [{\J}(X_k:B[X_1,\dots {\hat X}_k,\dots ,X_n]),X_k]$\\
(if the right-hand side exists).
 \item[j.7.] \ $j(A:B)=0\Leftrightarrow A,B$ are 
freely independent.
\end{description}

The {\it liberation Fisher information} $\var^*$ is defined by
\[
\var^*(A:B) \ = \ |j(A:B)|^2_2
\]
if $j(A:B)$ exists and $=+\i$ otherwise.

Among its properties is an inequality for freely Markovian triples $A,B,C$
which resembles the Stam inequality
\[
\var^*(A:C)^{-1} \ \geq \ \var^*(A:B)^{-1}+\var^*(B:C)^{-1} \ .
\]
Finally, the {\it mutual free information} $i^*$ is then given by
\[
i^*(A:B) \ = \ 
\textstyle{\frac 12} \int^\i_0 \var^*(U(t)AU(t)^*: B)dt
\]
where $(U(t))_{t\geq 0}$ is the unitary free Brownian motion which is
free w.r.t.~$A\vee B$.

\subsection{A variational problem for $\chi(\xr)$}

It is a natural variational problem for the free entropy to maximize
$$
\chi(\xr )-\t(P(\xr))  \eqno(*)
$$
where $X_j=X^*_j\in (M,\t)$, $\lejn$ and $P=P^*\in{\C}\<t_1,\dots ,t_n\>$
(see $\chi.7$ in 2.3 for this notation). The question is to find the joint
distribution of $(\xr)$ for which $(*)$ is maximum.
($(M,\t)$ is a ``universal" II${}_1$ factor containing all separable
II${}_1$ factors.)

It is interesting to note that this problem, about which we know very little
in this generality, appears to be connected to an important class of
\rn matrix models, about which similarly   very little is known
in full generality. To explain this, we shall consider the critical point
condition, which is a consequence of $(\xr)$ being a point where the maximum
is attained:
\begin{eqnarray*}
&&\frac{d}{d\e} \ \chi(X_1+\e P_1(\xr),\dots ,X_n+\e P_n(\xr))|_{\e =0}\\
&& \qquad =
\frac{d}{d\e} \ \t(P(X_1 +\e P_1(\xr),\dots ,X_n+\e P_n(\xr))|_{\e =0}
\end{eqnarray*}
Let $\p_j$ denote $\p_{X_j:{\C}[X_1,\dots {\hat X}_j,\dots X_n]}$
and let $d_j$ denote the cyclic derivative w.r.t.~$X_j$, i.e.,
\[
d_j \ = \ m\ \circ \ \sim \ \circ  \ \p_j 
\]
where $\sim$ is the flip for  ${\C}[\xr ]\ot{\C}[\xr]$ and
\[
m:{\C}[\xr ]\ot {\C}[\xr ]\to {\C}[\xr ] 
\]
is multiplication.
Then the critical point condition in view of $\chi.7$ becomes
\[
\sum_{\lejn} (\t\ot\t)(\p_j P_j) \ =\sum_{\lejn} \t((d_jP)P_j) 
\]
which in view of 3.3 means precisely that the conjugate variables
\[
{\J}_k \ = \ {\J}(X_k: {\C}[X_1,\dots {\hat X}_k,\dots X_n]
\]
exist and that
$$
{\J}_k \ = \ d_kP \qquad 1\leq k\leq n  \ . \eqno(**)
$$
Note that an equivalent way of stating these conditions is
$$
\sum_{i_j=k}\t(X_{i_1}\dots X_{i_{j-1}})\t(X_{i_{j+1}}\dots X_{i_p})
=\t(X_{i_1}\dots X_{i_p}(d_kP)(\xr))  \eqno(***)
$$
for all $1\leq k\leq n$ and monomials $X_{i_1}\dots X_{i_p}$.

The same equations (see [9] for instance) appear in the study 
of the large $N$ limit
of the general \rn multi-matrix model arising from a \pr ty measure with
density
\[
c_N \ e^{-N{\mb{Tr}}\, P(A_1,\dots ,A_n)}
\]
on the space of $n$-tuples of hermitian $N\x N$ matrices.

Like in the study of \rn matrix models also for the variational problem $(*)$,
it is natural to assume certain lower bounds for $P$. For instance the
condition
\[
\t(P(\xr))\geq A+B\log \t(X^2_1+\dots + X^2_n)
\]
where $B>\frac 12$ combined with $\chi.1$ gives
\[
\chi(\xr)-\t(P(\xr))\leq K-\e \log\t(X^2_1+\dots + X^2_n)
\]
for some constants $K$ and $\e>0$ which then will give a bound on $\t
(X^2_1+\dots + X^2_n)$ for a maximum.

For the reader familiar with one \rn matrix models, let us point out that for
$n=1$, the variational problem $(*)$ with $\mu =\mu_X$ the
distribution of $X=X_1$, becomes in view of $\chi.4$:
\[
\iint \log |s-t|d\mu(s)d\mu(t) -\int P(t)d\mu(t)
\]
while the equation $(**)$ in view of ${\J}.1$ becomes
\[
(2\pi H\mu)(X) \ = \ P'(X)
\]
or equivalently
\[
2\pi H\mu(t) \ = \ P'(X)
\]
$\mu$ -- almost everywhere (i.e., under continuity conditions for 
$t\in{\mb{supp}}\,\mu)$. 

These are familiar objects in the study of 1-\rn matrix models in the large
$N$ limit and free \en appears to provide the generalization of these for
multi-matrix models.

\subsection{Comments}

In this section we briefly discuss some of the problems encountered in the
effort to complete the theory and we also briefly mention further work in this
area, not covered in the previous sections.

{\bf Unification problem.} The ultimate goal of a complete theory also
would mean unification of the matricial \mi approach, the infinitesimal
approach and the mutual free information defined using $\d_{A:B}$
and the liberation process. This would mean in particular proving general
results of the form
\[
\chi(\xr) =\chi^*(\xr)
\]
and 
\begin{eqnarray*}
&& i^*(W^*(\xr):W^*(Y_1,\dots ,Y_m)) \\
&& \qquad = \chi(\xr)+\chi(Y_1,\dots ,Y_m) -\chi(\xr,Y_1,\dots ,Y_m)
\end{eqnarray*}
(when the $\chi$'s in the right-hand side are finite).

Clearly such results are a long way to go from where the theory is at present.
As always skeptics would raise the perspective of a negative answer.
On the other hand the results paralleling the classical theory, obtained thus
far, coupled with our general faith in beautiful mathematical theories should
be reasons for optimism that some form of a complete theory and unification
are possible. From a more pedestrian point of view it is clear that
unification will also very much depend on solving the technical problems in
completing each of the three directions.

\bs{\bf Technical problems.} Developing free \en theory in the  infinitesimal
approach, the problems one is facing at present are ``free analysis"
questions. Here is perhaps the simplest continuity question one would like to
settle in the affirmative:
\begin{quote}
{\it is \ ${\J}(X_1+tS_1,\dots ,X_n+tS_n) \in L^2(M,\t)$ 
a continuous function of} $t\in (0,\i)$?
{\it Here $S_1,\dots ,S_n$ are $(0,1)$-semicircular on
$\{S_1\},\dots ,\{S_n\},\{\xr\}$ are freely independent in}
$(M,\t)$.)
\end{quote}
The question is equivalent to the apparently weaker question: 
\begin{quote}
{\it is \ $\Phi^*(X_1+tS_1,\dots ,X_n+tS_n)$ 
as a  function of $t\in (0,\i)$ right continuous?
It is known the function is left continuous and decreasing.}
\end{quote}
Under this form the one-variable case, $n=1$, has been answered in
the affirmative in [41].

The problem of establishing a change of variables results for
$\chi^*(\xr)$ also runs into difficulties, part of which are related to
continuity questions like the preceding one.

\bs{\bf Free Fisher information relative to a completely positive map.}
Several results in the infinitesimal approach have been shown to hold in a 
more general framework involving a unital completely positive map
$\eta:B\to B$ ([25]). Instead of letting $\p_{X:B}$ take values in $M\ot M$
endowed with the scalar product derived from $\t\ot\t$ one uses the scalar
product
\[
\<x_1\ot x_2,y_1\ot y_2\> \ = \ \t(x^*_2\eta (E_B(x^*_1y_1))y_2) \ .
\]
One context where this generalization has a natural \mi counterpart occurs in
the study of Gaussian \rn band matrices [24],[17]. Another context involves
measure-preserving equivalence relations, and a free \pr ty interpretation
[26] of the recent work on the cost of such equivalence relations [14].

\bs{\bf Large deviations.}
Recent work on large deviations of Gaussian \rn matrices, up to technical
differences, can be viewed as aiming to prove a strengthening of the equality
of the free \en via \mi $\chi$ to the free \en $\chi^*$ defined via an
infinitesimal approach, i.e., a strengthening of the unification problem.
Slightly more precisely, the asymptotic of $k^{-2}\log \
{\mb{vol}} \ \G(k)$ where $\G(k)$ is a set of matricial \mi specified by
giving intervals for a finite number of normalized \ncc moments, should be
evaluated by the supremum of a rate function, involving the free \en $\chi^*$,
over the $n$-tuples of hermitian operators in tracial $W^*$-\pr ty spaces
satisfying the moment conditions. Even more precisely, the preceding should be
amended by taking Gaussian measure, removing cutoffs, replacing usual moments
by traces of products of some \ncc resolvents, etc.

In the one-variable case, both free \en [35 II] and the 
 large deviation question [1] are completely clarified and fit quite well
together.  In several variables a complete large deviations result, up to some
technical differences on microstates, would imply affirmative answers to the
lim sup, versus lim inf problem in Remark a) of 2.2 and of the Connes problem
in Remark d) of 2.2. Having in mind that a full  large deviations would imply
the solution of these difficult problems, note that the $n$-variable results
in [5] provide at present the closest result to a majorization of $\chi$ by
$\chi^*$. Besides the technical differences concerning microstates pointed out
above, there is one more important modification in [5] to be pointed out:
$\Phi^*$ is modified by the $L^2$-distance of
$({\J}(X_k:{\C}[X_1,\dots ,{\hat X}_k,\dots ,X_n]))_{1\leq k\leq n}$
to the set of cyclic gradients. This leads naturally to the problem whether
this $L^2$-distance is zero, i.e., whether the modification of $\Phi^*$
 is not really a modification of the quantity? Very little is known about
this. A purely \alg ic result in [38] implies the distance is zero
when the partial free Brownian gradients
$({\J}(X_k:{\C}[X_1,\dots ,{\hat X}_k,\dots ,X_n])$ are \ncc polynomials in
$\xr$. In a forthcoming paper by T.~Cabanal-Duvillard and A.~Guionnet
it is shown that the $n$-tuples of \ncc random variables for which
 Connes' problem
has an affirmative answer, are in the closure in distribution of those
 for which the above
question has an affirmative answer.

In another direction it is important to note that the large deviation work [5]
has brought powerful stochastic analysis techniques, applied to
matricial Brownian motions, to bear on the problems in this area.

\bs{\bf Some extremal problems.} Important classes of operators in
II${}_1$-factors, like the circular elements, are the solution to 
extremal problems for \en [21].

\bs{\bf The co\alg \ of} $\p_{X:B}$. \ The derivation of $\p_{X:B}$
is a comultiplication for a co\alg \ structure on $B[X]$. This leads to a
class of co\alg s where the comultiplication is a derivation, which has
remarkable duality properties closely related to results on conjugate
variables ${\J}(X:B)$ ([31]).

\clearpage
\begin{center}{\bf References}\end{center}

\begin{description}\item[[\,1\!\!]]
B.Ben Arous and A.Guionnet, Large deviations for
Wigner's law and Voiculescu's \ncc entropy. {\it Prob. Th. Rel. Fields} {\bf
108} no.~4 (1997), 517--542.
\item[[\,2\!\!]] P.Biane, Free Brownian motion, free stochastic calculus
and \rn matrices, in [40], pp.1--19.
\item[[\,3\!\!]] P.Biane and R.Speicher, Free diffusions, free \en and free
Fisher information, preprint (1999).
\item[[\,4\!\!]] P.Biane and D.Voiculescu, A free \pr ty analogue of the
Wasserstein metric on trace-state space, preprint.
\item[[\,5\!\!]] T.Cabanal-Duvillard and A.Guionnet, Large deviations,
upper bounds and \ncc entropies for some matrices ensembles, preprint.
\item[[\,6\!\!]] T.M.Cover and J.A.Thomas, {\it Elements of Information
Theory}. Wiley Interscience Publishers (1991).
\item[[\,7\!\!]] J.Dixmier, {\it Les $C^*$-alg\`ebres et leurs
Repr\'esentations.} Gauthier-Villar, Paris (1964).
\item[[\,8\!\!]] J.Dixmier, {\it Les Alg\`ebres d'Operateurs dans l'Espace
Hilbertien.} Gauthier-Villar, Paris (1969).
\item[[\,9\!\!]] M.Douglas, Large $N$ quantum field theory and matrix models,
 in [40], pp. 21--40.
\item[[\,10\!\!]] K.J.Dykema, Free products of hyperfinite von Neumann \alg s
and free dimension. {\it Duke Math. J.} {\bf 69} (1993), 97--119.
\item[[\,11\!\!]] K.J.Dykema, On certain free product factors via an extended
matrix model. {\it J. Funct. Anal.} {\bf 112}, 31--60.
\item[[\,12\!\!]] K.J.Dykema, Two applications of free entropy. {\it Math.
Ann.} {\bf 308} (1997), 547--558.
\item[[\,13\!\!]] J.Feldman and C.C.Moore, Ergodic equivalence relations,
cohomology and von Neumann \alg s, I, II. {\it Trans. Amer. Math. Soc.} {\bf
234} (1977), 289--359.
\item[[\,14\!\!]] D.Gaboriau, Co\^{u}t des relations d'equivalence et des
groupes. {\it Invent. Math.} {\bf 139} (2000), 41--98.
\item[[\,15\!\!]] L.Ge, Applications of free \en to finite von Neumann \alg s.
{\it Amer. J. Math.} {\bf 119} (1997), 467--485; ibidem II,
{\it Ann. of Math.} {\bf 147} (1998), 143--157; with J.Shen, ibidem III,
preprint.
\item[[\,16\!\!]] L.Ge and S.Popa, On some decomposition properties for
factors of type II${}_1$. {\it Duke Math. J.} {\bf 94} (1998), 79--101.
\item[[\,17\!\!]] A.Guionnet, Large deviations, upper bounds and central limit
theorems for band matrices and \ncc functionals of Gaussian large \rn
matrices, preprint.
\item[[\,18\!\!]] F.Hiai and D.Petz, Eigenvalues density of the Wishart matrix
and large deviations. {\it Infinite Dim. Anal. Quantum Prob.} {\bf 1} (1998),
633--646.
\item[[\,19\!\!]] R.Kadison and J.Ringrose, {\it Fundamentals of the Theory of
Operator Algebras} (3 volumes). Birkh\"auser, Boston.
\item[[\,20\!\!]] V.A.Marchenko and L.A.Pastur, The distribution of eigenvalues
in certain sets of \rn matrices. {\it Math. Sb.} {\bf 72} (1967), 507--536.
\item[[\,21\!\!]] A.Nica, D.Shlyakhtenko and R.Speicher, Some minimization
problems for the free analogue of the Fisher information.
{\it Adv. Math.} {\bf 121} (1999), 282--347.
\item[[\,22\!\!]] F.Radulescu, Random matrices, amalgamated free products and
subfactors of the von Neumann \alg \ of a free group, of noninteger index.
{\it Invent. Math.} {\bf 115} (1994), 347--389.
\item[[\,23\!\!]] C.E.Shannon and W.W.Weaver, {\it The Mathematical Theory of
Communication.} University of Illinois Press, Urbana, IL (1949).
\item[[\,24\!\!]] D.Shlyakhtenko, Random Gaussian band matrices and freeness
with amalgamation. {\it International Math. Res. Notices} no.~{\bf 20} (1996),
1013--1025.
\item[[\,25\!\!]] D.Shlyakhtenko, Free \en with respect to a completely
positive map. {\it Amer. J. Math.} {\bf 122} (2000), 45--81.
\item[[\,26\!\!]] D.Shlyakhtenko, Free  Fisher information with respect to a
completely positive map and cost of equivalence relations. MSRI preprint
1999--030 (1999).
 \item[[\,27\!\!]] D.Shlyakhtenko, On prime factors of type III.
{\it Proc. Nat. Acad. Sci.} {\bf 97} (2000), 12439--12441.
\item[[\,28\!\!]] R.Speicher, Combinatorial theory of the free product 
with amalgamation and operator-valued free \pr ty theory.
{\it Memoirs of the AMS} {\bf 627} (1998).
\item[[\,29\!\!]] S.Stratila and L.Zsido, {\it Lectures on von Neumann
Algebras.} Editura Academiei and Abacus Press (1979).
\item[[\,30\!\!]] M.B.Stefan, The indecomposability of free group factors over
nonprime subfactors and abelian sub\alg s, preprint.
\item[[\,31\!\!]] S.V.Szarek and D.Voiculescu, Volumes of restricted Minkowski
sums and the free analogue of the \en power inequality. {\it Comm. Math. Phys.}
{\bf 178} (1996), 563--570.
\item[[\,32\!\!]] D.Voiculescu, Symmetries of some reduced free product
$C^*$-\alg s, in {\it Operator Algebras and Their Connections with Topology
and Ergodic Theory}, Lecture Notes in Math., vol.~{\bf 1132}, Springer (1985),
pp.556--588.
\item[[\,33\!\!]] D.Voiculescu, Limit laws for \rn matrices and free products.
{\it Invent. Math.} {\bf 104} (1991), 201--220.
\item[[\,34\!\!]] D.Voiculescu, Circular and semicircular systems and free
product factors, in {\it Operator Algebras, Unitary Representations, Enveloping
Algebras and Invariant Theory}, Progress in Mathematics {\bf 92}, Birkh\"auser
(1990), pp.45--60.
\item[[\,35\!\!]] D.Voiculescu, The analogues of \en and of Fisher's
information measure in free \pr ty theory, {\it Comm. Math. Phys.} {\bf 155}
(1993), 71--92;
ibidem II, {\it Invent. Math.} {\bf 118} (1994), 411--440; 
ibidem III: The absence of Cartan sub\alg s, {\it Geom. Funct. Anal.} {\bf 6},
no.~1 (1996), 172--199; ibidem IV: Maximum \en and freeness, in [40],
pp.293--302; ibidem V: Noncommutative Hilbert transforms, {\it  Invent. Math.}
{\bf 132} (1998), 182--227; ibidem VI: Liberation and mutual free information,
{\it Advances in Math.} {\bf 146} (1999), 101--166.
\item[[\,36\!\!]] D.Voiculescu, A strengthened asymptotic freeness result for
\rn matrices with applications to free entropy. {\it International Math. Res.
Notices} {\bf 1} (1998), 41--64.
\item[[\,37\!\!]] D.Voiculescu, The coalgebra of the free difference quotient
and free
\pr ty theory,  {\it International Math. Res.
Notices} no.~2 (2000), 79--106.
\item[[\,38\!\!]] D.Voiculescu, A note on cyclic gradients, preprint.
\item[[\,39\!\!]] D.Voiculescu, Lectures on free \pr ty theory, in
{\it Lectures on Probability Theory and Statistics},
Ecole d'Et\'e de Probabilites de Saint-Flour XXVIII (1998),
Springer Lecture Notes in Math.~{\bf 1738}, pp. 280--349.
\item[[\,40\!\!]] D.Voiculescu, editor, {\it Free Probability Theory.}
Fields Institute Communications, (1997) vol.~{\bf 12}, AMS, Providence, RI.
\item[[\,41\!\!]] D.Voiculescu,  The derivative of order $\frac 12$ of a free
convolution by a semicircle distribution. {\it Indiana Univ. Math. J.}
{\bf 46}, no.~3 (1997), 697--703.
\item[[\,42\!\!]] D.Voiculescu, K.J.Dykema and A.Nica, {\it Free Random
Variables}, CRM Monograph Series (1992), vol.~1,  AMS, Providence, RI.
\item[[\,43\!\!]]  E.Wigner, On the distribution of the roots of certain
symmetric matrices. {\it Ann. Math.} {\bf 67}, (1958), 325--327.
\end{description}

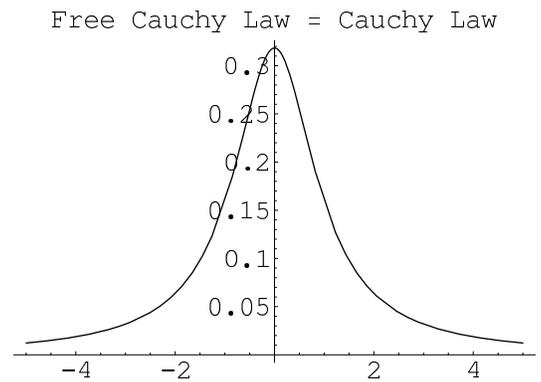
\end{document}